\documentclass[12pt]{amsart}
\usepackage{amssymb}
\usepackage{amsmath}
\usepackage{amsthm}
\usepackage{amsfonts}
\usepackage{stackrel}
\usepackage{enumitem}
\setlist{topsep=.3cm, itemsep=.2cm} 
\usepackage{mathrsfs}

\usepackage[T1]{fontenc}
\newtheorem{thm}{Theorem}[section]

\newtheorem{pro}[thm]{Proposition}
\newtheorem{cor}[thm]{Corollary}
\newtheorem{lem}[thm]{Lemma}

\newtheorem{example}[thm]{Example}
\numberwithin{equation}{section}
\usepackage{hyperref}
\newcommand{\R}{\mathbb{R}}

\newcommand{\D}{\displaystyle}

\def\qq#1{\qquad \mbox{#1}\quad}

\def\q#1{\quad \mbox{#1}\ }
\newcommand{\al}{\alpha}
\newcommand{\be}{\beta}

\newcommand{\de}{\delta}

\newcommand{\g}{\gamma}

\newcommand{\la}{\lambda}

\newcommand{\na}{\nabla}
\newcommand{\Om}{\Omega}

\newcommand{\Omb}{\overline{\Om}}
\newcommand{\p}{\partial}

\newcommand{\te}{\theta}

\title{Uniform  {\it a priori} bounds for Slightly Subcritical Elliptic Problems}
	
\author{Mabel Cuesta}
\address{Department of Mathematics, Université du Littoral Côte d'Opale (ULCO), Laboratoire de Mathématiques Pures et Appliqu\'ees Joseph Liouville (LMPA),		62100 Calais   (France)}
\email{\tt mabel.Cuesta@univ-littoral.fr}
\author{Rosa Pardo}
\address{Departamento de An\'alisis Matem\'atico y Matem\'atica Aplicada,  Universidad Complutense de Madrid, 28040{\textbf -}Madrid, Spain}
\email{\tt rpardo@ucm.es}

\thanks{Rosa Pardo is supported by grants  PID2022-137074NB-I00,  MICINN,  Spain, and by UCM, Spain,  Grupo 920894}

\begin{document}
\maketitle

\begin{abstract}
We obtain a uniform $L^{\infty}(\Omega)$ {\it a priori} bound, for  any positive weak solutions to   elliptic problem  with a  nonlinearity  $f$ {\it slightly subcritical, slightly superlinear}, and {\it regularly varying}.

To achieve our result, we first  obtain  a  uniform estimate of an specific $L^1(\Omega)$ weighted norm. This, combined with moving planes method and elliptic regularity theory, provides a uniform  $L^\infty$ bound  in a neighborhood of the boundary of $\Omega$.
Next, by using Pohozaev's identity, we obtain a uniform   estimate 
of one weighted norm of the solutions.
Joining now elliptic regularity  theory, and Morrey's Theorem,  we estimate from below the radius of a ball where a solution exceeds the half of its $L^\infty(\Omega)$-norm.
Finally, going  back to the previous uniform   weighted norm  estimate,
we conclude our result.


\medskip

\end{abstract}

\section{Introduction}
The aim of the present work is to  study  the existence of uniform $L^\infty(\Om)$ {\it a priori} bounds for positive solutions to  the following boundary-value problem:
\begin{equation}
\label{eq:ell:pb} 
\left\{ 
\begin{array}{rcll} 
-\Delta u&=&f(u), & \qquad \mbox{in } \Omega, \\ 
u&=& 0, & \qquad \mbox{on } \partial \Omega,
\end{array}\right.
\end{equation}
where $\Omega \subset \R ^N $ is a bounded  domain  with $C^{2,\al}$  boundary, with $N>2,$
and $f \in C\big(\R^+,\R^+\big)$   is locally Lipschitz, superlinear and of the form
\begin{align}
\label{defL}
&f(s)=s^{q}L(s),\quad s>0,\qq{with} q\in [1, 2^*-1],
\end{align}
where $L$ is a  positive, $C^1$-function for $s>0$, that is 
{\it slowly varying}   at infinity, (see \eqref{reva} for a precise definition with $q=0$). 
We prove a priori bounds for any {\it slightly superlinear} $f$ (see \ref{f2}) satisfying \eqref{defL} with $q\in[1,2^*-1)$. In case $q=2^*-1$, we prove uniform a priori bounds under a  condition involving $L'$ (see \ref{f3}).
\smallskip

Solutions to \eqref{eq:ell:pb} will be understood in the weak sense, although classical regularity results imply that weak solution to \eqref{eq:ell:pb} are at least class  $C^{2}(\Omb)$.
\smallskip

The problem on {\it a priori} bounds in the $L^\infty(\Om)$-norm of positive solutions   is a longstanding open problem, raised  by Gidas and Spruck in \cite{Gidas_Spruck_bd} as well as by de Figueiredo, Lions and Nussbaum in  \cite{deF_Lion_Nus_1982}. 
\smallskip

As a general goal,  we would like to obtain {\it a priori} bounds  estimates  for solutions  with functions $f$ that are  {\it slightly subcritical} at infinity 
$$
\lim_{s\to+\infty} \frac{f(s)}{s^{2^*-1}}=0,
$$
and whose behavior can be quite far from power-like. This led to consider  several natural classes of non-linearities, including nonlinearities with {\it regular variation at infinity},  (see \eqref{reva} for a precise definition).  Typical examples of nonlinearities $f$ with regular variation of index $q>0$ denoted as $f\in RV_q$ are 
$$
f(s)=s^q\log^\al  (K+s) \quad \hbox{ with } K>1,\hbox{ and } \al\in \R. 
$$

In a recent work \cite{Souplet_DCDS_23}  Souplet deals with  problems  with nonlinearities  of this type to obtain,  when $q\in(1,2^*-1)$,  uniform estimates of the infinity norm of positive solutions
(see   \cite[Theorem 3.1]{Souplet_DCDS_23}).
 He also  mentions 
`it seems unknown whether some universal estimates are true for solutions to 
problems with more general factors such as iterated logarithms, etc'.  Let us point out that in \cite[Theorem 1.1 and Theorem 1.2]{deF_Lion_Nus_1982}, de Figuei\-redo, Lions and Nussbaum provide universal estimates for $q\in(1,2^*-1)$  that include examples with logarithms factors as above (see Section \ref{sec:2} and equations \eqref{hyp:deFLN} and \eqref{hyp:deFLN2} for a discussion on their hypothesis). 
\smallskip

The question of the existence of uniform $L^\infty(\Om)$ norms   in the slightly subcritical case was first  addressed by one of  the authors in \cite[Theorem 1.1 and Corollary 2.2]{Castro_Pardo_RMC_2015} for the particular case 
$$
f(s)=\frac{s^{2^*-1}}{\log^\alpha(\mathrm{e}+s)}, \quad s\geq 0, \qq{and} \alpha>\frac{2}{N-2}\,.
$$
This nonlinearity $f$ is an example of a slightly subcritical, regularly varying function at infinity of order $2^*-1$. 
\smallskip

Our aim is to generalize the previous results  in the case 
where   $q=2^*-1$, $f$ is slightly subcritical and slightly superlinear at infinity, and
$L$ is a  positive, 
$C^1$-function for $s>0$, that is decreasing for $s$ large, and  {\it slowly varying}   at infinity.
\medskip

Let us remark here that, for slightly subcritical 
functions   $f\in RV_{2^* -1}$, the blow-up argument  fails. Indeed, the argument consist of, given a sequence of unbounded positive solutions $u_k$, we define a rescaling sequence
$$
v_k(y):=M_k^{-1} u_k(x_k+\de_k y), \quad y\in \Omega_k:=\frac{1}{\de_k}(\Omega-x_k),
$$
with $M_k=u_k(x_k)=\|u_k\|_\infty$, and  either $\de_k=\big(\frac{M_k}{f(M_k)}\big)^{1/2}$
or $\de_k=M_k^{-\frac{2}{N-2}}$.
Roughly speaking, in the convex case their limit $0<v\le v(0)=1$ solves either $-\Delta v=v^{2^*-1},$ or $-\Delta v=0$ in $\R^N$. 
In the first case,  Gidas-Ni-Nirenberg proved that the only classical positive solution 
of the elliptic equation in the whole $\R^N$ with  finite energy, is  
the so-called standard bubble
\begin{equation*}
U_{\delta}(x)= \delta^{N-2}\,(\de^2+|x|^2)^{-\frac{N-2}{2}},\qquad  
\delta=\sqrt{N(N-2)}\;,
\end{equation*}
see \cite{Gidas_Ni_Nirenberg}. In the second case, the Liouville Theorem, see \cite{Axler_Bourdon_Ramey, Nelson}, implies that $v\equiv 1$ in $\R^N$.
Consequently, the usual argument of contradiction based on non-existence fails.
\medskip 

Our approach to obtain a priori bounds for positive solutions,  combines  the Pohozaev's identity, the moving planes method,  the elliptic regularity, and Sobolev embeddings. 
\medskip

Precisely,   we will assume   the following hypothesis for $q\in [1,2^*-1]$ :

\begin{enumerate}[label=\textbf{(f\arabic*)$_\infty$}]
\item
\label{f1}
$f$  is {\it slightly subcritical} at infinity : ${\D\lim_{s\to +\infty}}\,\frac{f(s)}{s^{2^*-1}}=0.$

\item 
\label{f2}
$f$ is {\it slightly superlinear} at infinity in the sense of de Figueiredo, Lions and Nussbaum (see \cite[p. 43]{deF_Lion_Nus_1982})
\begin{equation*}
{\D\liminf_{s\to +\infty}}\, \frac{f(s)}{s}> \lambda_1,
\end{equation*}
where $\la_1$ is the first eigenvalue for the Dirichlet eigenvalue problem. 
\smallskip 

If the domain $\Om$ is not convex, we will assume moreover that the function $f(s)/s^{2^*-1}$ is nonincreasing for any $ s > 0$.\\

Besides, when   $q=2^*-1$, we will also assume :
\item
 \label{f3}  The function $L$ defined in \eqref{defL} is 
 $C^1(0,+\infty),$  and $L'$  satisfies:
\begin{align}
&L'(s)< 0,\quad \hbox{ for all }  s \ge s_1,  \hbox{ for some } s_1\gg 1;\label{L':neg}\\
&|L'|\in RV_{-1} \text{ that is : }{\D\lim_{s\to +\infty}}\, \frac{|L'(\tau s)|}{|L'(s)|} =\frac{1}{\tau}\qquad\forall  \tau >0.\nonumber
\end{align}
\end{enumerate}

A measurable  positive function $g:\R^+\to\R^+$ is said to be a {\it a regularly varying function (at infinity) of index $q\in\R$ } if 
\begin{equation}
\label{reva}
\lim_{s\to+\infty} \frac{g(\tau s)}{g(s)}=\tau^q \qquad  \forall \tau>0.
\end{equation}
In that case we write $g\in RV_q$. Positive functions satisfying condition \eqref{reva}  with $q=0$, are called {\it slowly varying (at infinity)}. The set of {\it slowly varying functions at infinity} is denoted as $RV_0$. 
We mention the seminal work of Karamata \cite{Karamata}, and the books  \cite{Bingham, Resnick, Seneta},
or for instance the following papers for several applications \cite{Costa_Quoirin_Tehrani, G-Melian_Iturriaga_RQuoirin}. In particular, by  the Uniform Convergence Theorem of Karamata, hypothesis \ref{f3} implies that $L\in RV_0$ (see  Theorem    \ref{usfF}.(i)).

\medskip

Obviously, asymptotic functions to  a non-trivial constant at infinity   are slowly varying functions.
But, a slowly varying function $L$ may exhibit oscillations with infinite amplitude:
$$
\liminf_{s\to\infty} L(s) = 0,\qquad \limsup_{s\to\infty} L(s) = +\infty,
$$
an example being 
$$
L(s)=\exp \Big[\sin \big(\log^\be (s)\big) \log^\be (s)\Big],\qq{with}0<\be<1/2.
$$ 
This kind of function is out of our framework, since it does not satisfy condition \eqref{L':neg}.
It is an open problem to find $L^\infty(\Om)$ {\it a priori} bounds without imposing either the monotonicity of $L$ at infinity,  \eqref{L':neg}, or,  for oscillatory nonlinearities, to be bounded from below by a monotone function.
The function 
$$
L(s)=\exp \Big[\al\Big(\g\sin \big(\log^\be (s)\big) +\log^\be (s)\Big)\Big],\text{ with }|\g|<1,\ 0<\be<1/2,
$$
is inside of our scope.
In the Appendix \ref{app:B} we show examples of nonlinearities $L$ satisfying \ref{f3} and hypothesis \eqref{hypothesis1} below. In particular, for oscillatory nonlinearities satisfying our hypothesis, see $L_{12}$ and $L_{13}$ defined in \eqref{12}-\eqref{13} respectively in Example \ref{example} of Appendix \ref{app:A}.\\

Our main result  is the following  Theorem: 

\begin{thm}{$L^\infty(\Om)$ uniform a priori bound for $q=2^*-1$.}
\label{th:apriori:cnys} 

Let $q=2^*-1$ and assume that $f:\R^+\to \R^+$  is a 
$C^1$ function   satisfying   hypothesis  {\rm \ref{f1}--\ref{f3}}.
If 
\begin{align}\label{hypothesis1}
& \lim_{s\to\infty} \,\frac{s|L'(s)|}{L^{N/2}(s)}=+\infty,
\end{align}
then, there exists a constant $C>0$ such that  for any  $u\in H_0^{1}(\Om)$ positive weak solution   to \eqref{eq:ell:pb},  
$$
\|u\|_{\infty}\le C, 
$$  
where  $C$ depends only on     $f$, $\Om$ and $N$,  and it is independent of the  solution $u$.
\end{thm}

We illustrate our approach by an exhaustive list of examples of regularly varying functions satisfying our hypothesis in Appendix \ref{app:A}  and Appendix \ref{app:B}.
Those examples are proposed in \cite[p. 47]{Seneta} and in \cite[p. 16]{Bingham}.
\begin{cor}
\label{cor:ex}
Let
\begin{align*}
&f_i(s)=s^{2^*-1}L_i(s), \q{for} s> 0 ,
\end{align*}
for  $L_i\in RV_0$,  with $i\in\{1,2,4-7,11-14\}$, depending on a parameter $\al<0,
$ and possibly more parameters, described in the Example {\rm \ref{example}} of the Appendix {\rm \ref{app:A}}. In case $i=1$, assume also that  $|\al|>\frac{2}{N-2}.$
\smallskip

Then, there exists a uniform constant $C>0$ such that for  every non-negative weak solution  $u$ of \eqref{eq:ell:pb}$_i$ (the problem \eqref{eq:ell:pb} for $f=f_i$), 
\begin{equation*}
\|u\|_{L^{\infty}(\Om)}  \leq C
\end{equation*}  
where $C$ depends only on   $f_i$,  $N$, and $\Om$, but it is independent of the  solution $u$.   
\end{cor}
See   the Appendices \ref{app:A} and \ref{app:B} for a proof.\\

We also present one alternative proof of a known result.

\begin{thm}
\label{th:apriori} {\rm $L^\infty(\Om)$ uniform a priori bound for $q\in [1,2^*-1)$.}

Let $q\in [1,2^*-1)$. Assume that $f:\R^+\to \R^+$  is a 
continuous function  satisfying   hypothesis  {\rm \ref{f1}--\ref{f2}}. Assume also that $L\in RV_0.$

Then, there exists a uniform constant $C>0$ such that for  every non-negative weak solution  $u$ of \eqref{eq:ell:pb}, 
\begin{equation*}
\|u\|_{L^{\infty}(\Om)}  \leq C
\end{equation*}  
where $C$ depends only on   $f$,  $N$, and $\Om$, but it is independent of the  solution $u$. 
\end{thm}

\medskip

In particular,  Theorem \ref{th:apriori} applies to the following nonlinearities. Let
\begin{align*}
&f_i(s)=s^{q}L_i(s), \ \text{ and either} 
\begin{cases}
q=1, &
\al>0,\\
q\in(1,2^*-1), &
\al\in\R,
\end{cases} 
\end{align*}
for $s> 0$, and for  $L_i\in RV_0$,  for $i=1,\cdots ,14$, depending on $\al,$ and possibly more parameters, described in the Example \ref{example} of the  Appendix \ref{app:A}. \\

 Let us   sketch our approach for the particular case when $L'<0$. First, using the moving planes method, as in \cite{deF_Lion_Nus_1982}, we  obtain uniform estimates in a neighborhood of the boundary for any positive non-negative solution to \eqref{eq:ell:pb}, see Theorem \ref{th:comp:02}.
\smallskip

Next, let us consider the Nemitsky operator associated to  the function  $2^* F(s)-sf(s)$, where $F(s):=\int_0^s f(t)dt$. On the one hand, Pohozaev's identity and the uniform estimates in a neighborhood of the boundary implies that $|\int_\Om 2^* F(u)-uf(u)\,dx| $ is uniformly bounded for any positive solution. On the other hand, in terms of $L$, one has for $s$ large
\begin{equation*}
2^* F(s)-sf(s)
=\int_0^s t^{2^* } |L'(t)|\,dt \sim  s^{2^*+1}|L'(s)|.    
\end{equation*}
Hence, we infer the uniform  estimate 
\begin{equation}\label{H1:bd}
\int_\Omega u^{2^*+1}|L'(u) |\,dx \leq C.    
\end{equation}
\smallskip

Finally, combining elliptic regularity  theory, with Morrey's Theorem,  and properties of slowly varying functions, see Corollary \ref{cor:RV:db}, we estimate from below the radius of a ball where a solution exceeds the half of its $L^\infty(\Om)$-norm, see \cite{Castro_Pardo_RMC_2015} for this technique.
Finally, going  back to the  previous bound \eqref{H1:bd},
we deduce that
\begin{equation*}
 \dfrac{\|u\|_\infty\ \big|L'\big(\|u\|_\infty\big)\big|}{L^\frac{N}{2}\big(\|u\|_\infty\big)}\le C<+\infty, 
\end{equation*}
which joint with hypothesis \eqref{hypothesis1}, implies that $\|u\|_\infty$ is bounded, and concluding our result, see Theorem \ref{th:apriori:cnys}.
\medskip

This paper is organized as follows. In Section \ref{sec:2} we review the result of \cite{deF_Lion_Nus_1982} to obtain a uniform $L^\infty(\Om)$ bound for positive solution,
and show why one of their  hypothesis, based on   Pohozaev's identity,  
can not be used in our case.  In Section \ref{sec:3} we  collect some known results. In particular, we include results of \cite{deF_Lion_Nus_1982}  proving a uniform $L^\infty$ bound  {\it near the boundary} of $\Omega$. Also in Section \ref{sec:3} we recall some properties of regularly  varying functions that we will use later.
In Section \ref{sec:4} we  prove Theorem \ref{th:apriori:cnys}, and in Section \ref{sec:5} we prove Theorem \ref{th:apriori}. We conclude this work with two Appendices. In Appendix \ref{app:A}  we provide examples of regularly varying  functions for which our theorems apply.  In Appendix \ref{app:B} we check all our hypothesis for the examples given in Appendix \ref{app:A}, in order to prove Corollary  \ref{cor:ex}.

\section{On de Figuei\-redo, Lions and Nussbaum condition}
\label{sec:2}
In \cite{deF_Lion_Nus_1982} de Figuei\-redo, Lions and Nussbaum firstly obtain a uniform $L^\infty$ {\it a priori} bound in a neighborhood of the boundary, combining a uniform $L^1(\Om)$ bound in the way of Brezis and Turner \cite{Brezis_Turner}, with the moving planes method of Gidas, Ni and Nirenberg \cite{Gidas_Ni_Nirenberg}. Secondly, 
using the Pohozaev identity, they obtain a uniform $H^1_0(\Om)$ {\it a priori} bound. Thirdly, using a device due to  Brezis and  Kato \cite{Brezis_Kato}, they obtain a uniform $L^q(\Om)$ {\it a priori} bound, for each $q<+\infty$. Finally, using elliptic regularity, they conclude the proof.\\

Let us review their hypothesis. In \cite{deF_Lion_Nus_1982}, the nonlinearity $f$ is assumed to satisfy
\begin{equation}\label{hyp:deFLN}
 \liminf_{s \to +\infty}\, \frac{ \theta F(s)-sf(s)}{s^2f(s)^{2/N}} \geq 0,\qquad\mbox{for some}\quad\theta \in [0,2^\star ),   
\end{equation}
where $F(s)=\int_0^s f(t)\, dt.$ They conjecture that  this condition is not necessary, but it is essential in proving their results.
It can be seen that for any
\begin{equation}\label{f:q:L}
f(s)=s^qL(s) \qq{with} L\in RV_0,        
\end{equation}
Theorem \ref{usfF} in section 3 implies that
\begin{equation}\label{F:q:te}
\lim_{s \to +\infty} \frac{ \theta F(s)-sf(s)}{sf(s)} =\frac{\te}{q+1}-1\ge 0\iff  q\le\te-1<2^\star-1.
\end{equation}
Moreover, 
$$
\liminf_{s \to +\infty} \frac{ sf(s)}{s^2f(s)^{\frac{2}{N}}} =\liminf_{s \to +\infty}s^{q\frac{N-2}{N}-1}L^{\frac{N-2}{N}}(s)=
\begin{cases}
+\infty &\text{if } q>\frac{N}{N-2}\\
\D\liminf_{s \to +\infty}L^{\frac{N-2}{N}}(s) &\text{if } q=\frac{N}{N-2}\\
0 &\text{if } q<\frac{N}{N-2},
\end{cases}
$$
see Proposition \ref{pro:L:al}.(ii).

Hence, hypothesis \eqref{hyp:deFLN} is satisfied for $q<2^*-1$, in that case
\begin{equation}\label{hyp:deFLN2}
\liminf_{s \to +\infty} \frac{ \theta F_i(s)-sf_i(s)}{s^2f(s)^{2/N}} =
\begin{cases}
+\infty &\text{if } q\in\left(\frac{N}{N-2},2^*-1\right)\\
\D\liminf_{s \to +\infty}L^{\frac{N-2}{N}}(s)\ge 0 &\text{if } q=\frac{N}{N-2}\\
0 &\text{if } q<\frac{N}{N-2}.
\end{cases}   
\end{equation}
\\

Assume now that $q=2^*-1$ for some $f$ given by \eqref{f:q:L}. For any $\te\in[0,2^\star ),$ multiplying and dividing by $sf(s)$, \eqref{F:q:te} and Proposition \ref{pro:L:al} imply that
\begin{align*}
\liminf_{s\to+\infty}\frac{\te F(s)-sf(s)}{s^2f(s)^{2/N}}&=\left(\frac{\te}{2^*}-1\right)\liminf_{s\to+\infty}\frac{sf(s)}{s^2f(s)^{2/N}}\\
& =\left(\frac{\te}{2^*}-1\right)\liminf_{s\to+\infty} s^{2/N}L(s)^{(N-2)/N}=-\infty.
\end{align*}
Hence, hypothesis \eqref{hyp:deFLN} is not satisfied,
and their method do not cover the slightly subcritical case with $q=2^*-1$.

\section{Known  results}
\label{sec:3}
\subsection{Uniform estimates on the boundary }

We present  the following theorem on the existence of a uniform $L^\infty$ bound in a neighborhood of the boundary.

\begin{thm}\label{th:comp:02} Assume that $\Omega \subset \R ^N $ is a bounded domain  with $C^{2}$ boundary.  Assume that the nonlinearity $f$ satisfies {\rm \ref{f2}}.

If $u\in C^2(\overline{\Om})$ satisfies  \eqref{eq:ell:pb} and $u>0$ in $\Om$, then  there exists a constant $d_0>0$ depending only on $\Om$ and not on $L$, $q$ or $u$, and a constant $C$ depending only on $\Om$, $L$ and  $q$ but not on $u$, such that
\begin{equation*}
\max_{\Om\setminus \Om_{d_0}} u\leq C
\end{equation*}
where $\Om_{d_0}:=\{x\in\Om\ : \ d(x,\p\Om)>d_0\}.$
\end{thm}

\noindent {\bf Sketch  of the proof of Theorem \ref{th:comp:02}.}
If $\Omega$ is convex, we observe that, reasoning as in \cite{deF_Lion_Nus_1982},  any positive solution $u$  is locally increasing  in a fixed neighborhood of the boundary, following directions close to the normal direction. This provides $L^\infty$ bounds locally in a neighborhood of the boundary. See for instance \cite[Theorem A.1]{Castro_Pardo_RMC_2015}.
\medskip

If $\Omega$ is not convex,   we had assumed also that the function $f(s)/s^{2^*-1}$ is nonincreasing for any $ s > 0$. At this moment, we observe that, reasoning as in \cite[Theorem 1.2]{deF_Lion_Nus_1982},  the Kelvin transform of $u$ at  $x_0\in\p\Om$ is locally increasing   in the maximal cap of the transformed domain. This provides $L^\infty$ bounds for the Kelvin transform locally.
By a compactification process, we then translate this into $L^\infty$ bounds in a neighborhood of the boundary for any solution to the elliptic equation. See for instance \cite[Theorem B.3]{Castro_Pardo_RMC_2015} for more details.
\qed

\subsection{Some properties of regularly varying functions}

Trivially,  for  any function $g\in RV_q$ there exists $L\in RV_0$ such that   $g(s)=s^q L(s)$ for $s>0$ large. 
\medskip

The following result was proved by Karamata \cite{Karamata} in the continuous case and extended by Korevaar, van Aardenne-Ehrenfest and de Bruijn in the measurable case, see \cite{Korevaar_Aardenne-Ehrenfest_Bruijn}, for the following result see 
\cite[Theorem 1.1]{Seneta}  and \cite[Theorem 1.5.2]{Bingham}.

\begin{thm}[Uniform Convergence Theorem]
    \label{usfF}
Let $g\in C(0,+\infty)$ positive at infinity and $g\in RV_q$ for some $q\in\R$.  Let $G(s)=:\int_0^s g(t)\,dt$. 
 Then: 
\begin{enumerate}[label=\rm{(\roman*)}]
\item  $G$ is a function regularly varying at infinity of index $q+1$,
 and satisfies
\begin{align} 
\lim_{s\to+\infty}\frac{sg(s)}{G(s)}=q+1.
\label{eq:usfF}
\end{align}
\item  Furthermore, for all $a,  b\in \R^+,$ $0<a<b$,  the limit  \eqref{reva}  
is uniform for $ \tau \in [a,b].$
\item  
Moreover, if $g$ is bounded close to 0 and $q>0$,  then
the limit  \eqref{reva} is uniform for $ \tau \in (0,b]$, for any $b>0$. 
\end{enumerate} 
\end{thm}

The following result can be read in \cite[Theorem 2]{Lamperti_1958a}. 
\begin{pro}\label{pro:lamp}
Let $g$ be a positive function and continuously differentiable in $(s_1, \infty)$ for some $s_1>0$. If
\begin{equation*}
sg'(s)/g(s) \to q\in\R \qq{as}   s\to\infty,
\end{equation*}
then $g\in RV_q$.
\end{pro}

For the following result,  see \cite[Proposition 1.3.6]{Bingham} for parts (i)-(iv). The proof of (v) is trivial.

\begin{pro}\label{pro:L:al}
Let $L\in C(0,+\infty)$ positive and varying slowly at infinity. 

\indent {\rm (i) }
$L(s)^\al \in RV_0$ for every $\al\in\R.$

{\rm (ii)}
If  $\al>0,$
$\ s^\al L(s)\to +\infty,\ s^{-\al} L(s)\to 0\ (s\to +\infty).$

{\rm (iii)}
If $L_1,\ L_2$ vary slowly at infinity,  so do $L_1+L_2$, $L_1 L_2$, and, if $L_2(s)>0 $ for all $s>0$ big enough, also $L_1/L_2$. 

{\rm (iv)} Let $L_1,\ L_2$ vary slowly at infinity. If $L_2(s)\to+\infty $ as $s\to+\infty $, $L_1\circ L_2$  is slowly varying.

{\rm (v)} If for all $ \tau>0$, $L(\tau s)-L(s)\to 0 $ as $s\to+\infty $, $e^{L}$ is slowly varying. 
\end{pro}

Since  Theorem \ref{usfF}.(iii), and Proposition \ref{pro:L:al}, we deduce the following Corollary. 

\begin{cor}
\label{cor:RV:db}
Let    $g\in C([0,+\infty))\cap RV_q$, with $q\in\R$. 
Then there  exists a constant $C_1$ such that 
\begin{equation}\label{bef5}
C_1\le   \liminf_{s\to+\infty}\,  \dfrac{\min_{[s/2,s]}\, g}{g (s)}.
\end{equation}
Moreover, if $q>0$ and $g$ is bounded close to 0, then there exists $C_2>0$ such that 
\begin{equation} 
\label{bef6}
\limsup_{s\to+\infty}\,  \dfrac{\max_{[0,s]}\, g}{ g(s)}\le C_2.
\end{equation}
\end{cor}

\begin{proof}
Since $g\in RV_q$  and Theorem \ref{usfF}.(ii), there exists $s_0>0$ such that 
$$
\frac12 \tau^q g(s)\le g(\tau s)\le 2 \tau^q g( s),\quad \forall s\ge s_0,\ \forall\tau\in [1/2,1].
$$
Now, denote by $t=\tau s$, then 
$$
\left(\frac12\right)^{q+1} g(s)\le g(t)\le  2 g( s),\quad \forall t\in [s_0/2,s],\quad \forall s\ge s_0,
$$
and so,   inequality \eqref{bef5} holds for $C_1:=\left(\frac12\right)^{q+1}.$
\medskip

Assume now that  $q>0$ and  $g$ bounded close to 0. Since  Theorem \ref{usfF}.(iii), there exists $s_0>0$ such that 
$$
g(\tau s)\le  (\tau^q+1) g( s),\quad \forall s\ge s_0,\ \forall\tau\in (0,1].
$$
Now, denote by $t=\tau s$, then 
$$
g(t)\le  2 g( s),\quad \forall t\in (s_0,s],\quad \forall s>s_0,
$$
and so,  the second inequality holds for $C_2:=2.$
\end{proof}

\section{Proof of Theorem \ref{th:apriori:cnys}}
\label{sec:4}
We recall that Theorem \ref{th:apriori:cnys} is for  $q=2^*-1$, and in that case, by \ref{f1}, $L(s)\to 0$ as $s\to\infty$.  From now on, all throughout this paper, $C$ denotes several constants
independent of $u$.
\medskip

\begin{lem}
\label{lem:3.1}
Assume that $f$ satisfies   {\rm \ref{f1}--\ref{f3}}. 
Let $L$  be defined by \eqref{defL}.
\medskip

Then,  there exists a constant $C>0$ such that     for any  $u\in H_0^{1}(\Om)$  positive  weak solution   to \eqref{eq:ell:pb},
\begin{equation}\label{u:fu:C}
\int_\Om u^{2^*+1} |L'(u)|\, \,dx \le C.
\end{equation}
 \end{lem}
\bigskip

\noindent {\bf Proof of Lemma \ref{lem:3.1}.} We split the proof into 4 steps. First, we look for a boundary uniform estimate, next for a Pohozaev's type estimate, then an equivalence at infinity, and finally for the uniform estimate \eqref{u:fu:C}.\\
 
\noindent {\bf Step 1}. {\it Boundary uniform estimates}.
\smallskip

From Theorem \ref{th:comp:02} and de Giorgi-Nash type Theorems, see  \cite[Theorem 14.1]{L-U}
$$
\|u\|_{C^{0,\nu}(\Omega_{d_0/8}\setminus \Omega_{7d_0/8})} \leq C,\qquad \mbox{for any}\quad\nu\in (0,1),
$$
where $\Om_{t}:=\{x\in\Om\ : \ d(x,\p\Om)>t\}.$
Hence, from Schauder interior estimates, see \cite[Theorem 6.2]{G-T}, we get 
$$
\|u\|_{C^{1,\nu}(\Omega_{d_0/4}\setminus \Omega_{3d_0/4})} \leq C.
$$
Finally, combining $L^p$ estimates with Schauder boundary estimates, see \cite{Brezis_2011}, \cite{G-T}
$$
\|u\|_{W^{2,p}(\Omega\setminus \Omega_{d_0/2})} \leq C,\qquad\mbox{for any}\quad p\in(1,\infty).
$$
Consequently, there exists two constants $C,\delta_0>0$ independent of $u$ such that
\begin{equation}\label{bdd:3}
\|u\|_{C^{1,\nu}(\Omb\setminus \Om_{\delta_0})} \leq C,\qquad\mbox{for any}\quad\nu\in (0,1).
\end{equation}
\medskip

\noindent{\bf Step 2}. {\it  A Pohozaev's type estimate}.
\smallskip

 From a slight modification of Pohozaev identity, see \cite[Lemma 1.1]{deF_Lion_Nus_1982} and \cite{Pohozaev}, if $y\in \R^N $ is a fixed vector, then any positive solution $u$ of \eqref{eq:ell:pb} satisfies
\begin{equation}\label{F:f:3}
2^* \int_\Om F(u)\, \,dx -\int_\Om uf(u)\, \,dx=
\frac1{N-2}\int_{\p\Om} (x-y)\cdot\nu(x)\, |\nabla u|^2 \, dS,
\end{equation}
where $\nu =\nu(x)$ denotes the unit  outward  normal to $\partial \Omega$ at $x$. Since \eqref{bdd:3}, the RHS is uniformly bounded, so the LHS is uniformly bounded,
\begin{equation}\label{F:f:4}
\left| 2^* \int_\Om F(u)\, \,dx -\int_\Om uf(u)\, \,dx\right| \le C.
\end{equation}
\medskip

\noindent{\bf Step 3}. {\it An equivalent function to $2^*F(s)-sf(s)$ at infinity.} 
\smallskip

We claim that 
\begin{equation}
\label{reg:var:Poho:q=}   
\lim_{s\to+\infty}\frac{2^*F(s)-sf(s)}{s^{2^*+1}\big|L'(s)\big|}=\frac{1}{2^*}.
\end{equation}
Indeed,
integrating by parts
$$
F(s)=\int_0^s t^{2^*-1}L(t)\,dt
=\frac{1}{2^*}s^{2^*}L(s)-\frac{1}{2^*}\int_0^s t^{2^*}L'(t)\,dt.
$$
Consequently
\begin{align}
\label{poho:crit}   
\lim_{s\to+\infty}\frac{2^*F(s)-sf(s)}{s^{2^*+1}\big|L'(s)\big|}
=\lim_{s\to+\infty}-\,\frac{\D\int_0^s t^{2^*}L'(t)\,dt}{s^{2^*+1}\big|L'(s)\big|}.
\end{align}
On the other hand, let us define 
\begin{equation*}
g(s):=s^{2^*}|L'(s)|,\qquad G(s)=\int_0^s g(t)\,dt.
\end{equation*}

Since $|L'|\in RV_{-1}$,  by  definition of a regularly varying function, see \eqref{reva}, it follows that $g\in RV_{2^*-1}$ and finally, from  Theorem \ref{usfF}(i)   we have that $G\in RV_{2^*}$ and  
$$
\lim_{s\to+\infty}\frac{G(s)}{sg(s)}=\frac1{2^*}.
$$
Hence
\begin{align}
\label{eq:sgG:2} 
\lim_{s\to+\infty}\frac{\D\int_0^s t^{2^*}|L'(t)|\,dt}{s^{2^*+1}|L'(s)|}=\frac{1}{2^*}.
\end{align}
Using  finally  that  $L'(s)< 0$   for $s\ge s_1$, see {\rm \ref{f3}},
 the claim \eqref{reg:var:Poho:q=}  follows from \eqref{poho:crit} and \eqref{eq:sgG:2}.

\bigskip

\noindent{\bf Step 4}.  {\it A uniform estimate of $\int_\Om u^{2^*+1} |L'(u)| \,dx$}. 
\smallskip

From \eqref{reg:var:Poho:q=} we infer that there exists $s_0>0$ such that
\begin{equation*}
2^*F(s)-sf(s)\ge \frac12\,\frac1{2^* }\,   s^{2^*+1 }|L'(s)|,\quad  \forall  s>s_0.
\end{equation*}

Applying the above inequality to any positive solution, and integrating on $\Om$ we obtain that
$$
2^*\int_\Om F(u)\, \,dx -\int_\Om uf(u)\, \,dx
\ge \frac12\,\frac1{2^* }\, \int_\Om u^{2^*+1}|L'(u)|\, \,dx-C_0,
$$
for some constant $C_0$ independent on $u$.  This estimate, 
the uniform $C^{1,\nu}$ estimate in a neighborhood of the boundary, see
\eqref{bdd:3}, 
and Pohozaev identity \eqref{F:f:3}, yield \eqref{u:fu:C}
for some constant $C$ independent of $u$. 
\qed
\bigskip

For the proof of Theorem \ref{th:apriori:cnys}, we estimate from below the radius of a ball where a solution exceeds the half of its $L^\infty(\Om)$-norm, see \cite{Castro_Pardo_RMC_2015} for this technique. To do that, we use elliptic regularity  theory $W^{2,q}(\Om)$ for $q>N/2$,  Sobolev embeddings,
Morrey's Theorem,  and properties of slowly varying functions. Going  back to the \eqref{u:fu:C}-estimate in Lemma \ref{lem:3.1},  we conclude that, for an unbounded sequence $\{u_n\}_n \subset H_0^{1}(\Om)$ of solutions  to \eqref{eq:ell:pb}, the following holds
$$
\dfrac{\|u\|_\infty\ \big|L'\big(\|u\|_\infty\big)\big|}{L^\frac{N}{2}\big(\|u_n\|_\infty\big)}\le C<+\infty,
$$ 
contradicting the hypothesis \eqref{hypothesis1}, and implying  our result.

\begin{proof}[Proof of Theorem \ref{th:apriori:cnys}]
Let  $\{u_n\}_n \subset H_0^{1}(\Om)$ be an unbounded sequence of solutions to \eqref{eq:ell:pb}, i.e.  $\|u_n\|_\infty\to +\infty$.  We can assume that  
$$
\|u_n\|_\infty \ge 2 s_1, 
$$
where $s_1>0$ comes from {\rm \ref{f3}}. 

By hypothesis $L'(s)<0$ for $s>s_1,$  so $|L'(u_n)|^{-1}$ is well defined for $u_n(x)\ge s_1$ and 
$$
\int_{\Om}\big|f(u_n)\big|^q\, \,dx=\int_{\Om\cap u_n(x)> s_1}\big|f(u_n)\big|^q\, \,dx+\int_{u_n(x)\le s_1}\big|f(u_n)\big|^q\, \,dx.
$$ 
Moreover,
\begin{align*}
&\int_{\Om\cap u_n(x)> s_1} \big|f(u_n)\big|^q\, \,dx\\
&\qquad =\int_{\Om\cap u_n(x)> s_1} \Big(u_n^{(2^*-1)q-(2^*+1)}L(u_n)^q |L'(u_n)|^{-1}\Big) \ u_n^{2^*+1}|  L'(u_n)|
dx  
\end{align*}
Let us denote 
$$
h(s):= s^{(2^*-1)q-(2^*+1)}L(s)^q |L'(s)|^{-1},\quad s\ge s_1.
$$
Hence $h\in RV_{2^*(q-1)-q}$,   with $2^*(q-1)-q>0\Longleftrightarrow q>2N/(N+2)$, which holds whenever $q>N/2$.  
Then by \eqref{bef6} in Corollary \ref{cor:RV:db},  and Lemma \ref{lem:3.1}
$$
\int_{\Om} \big|f(u_n)\big|^q\, \,dx
\le C \big(1+h\big(\|u_n\|_\infty\big)\big)\le Ch\big(\|u_n\|_\infty\big).
$$

Now, by Sobolev embeddings and elliptic regularity
\begin{align*}
\|u_n\|_{W^{1,q^*}(\Om)}&\le C\|u_n\|_{W^{2,q}(\Om)}\le C\left(\int_{\Om} \big|f(u_n)\big|^q\right)^\frac{1}{q}
\le  C h\big(\|u_n\|_\infty\big)^{1/q},\nonumber 
\end{align*}
where $\frac{1}{q^*}=\frac{1}{q}-\frac{1}{N}.$ Observe that $q^*>N.$\\

From Morrey’s Theorem, (see \cite[Theorem 9.12 and Corollary 9.14]{Brezis_2011}), there exists a constant $C$ only dependent on $\Om,$ $ q$ and $N$ such that
\begin{equation}\label{Morrey}
|u_n(x)-u_n(y)| \le C   \|u_n\|_{W^{1,q^*}(\Om)} \, |x-y|^{1-N/q^*},\qquad \forall x,y\in\Om.
\end{equation}
Next,  let us take  $x_n \in \Om$ such that
\begin{equation}\label{xn}
   u_n(x_n) =  \max_{\Om} u_n=\|u_n\|_\infty\ge  2s_1.
\end{equation}
Let  $B_n:=B(x_n,R_n)\subset\Om$, where $R_n,$ is defined such that
\begin{equation*}
u_n(x)\ge \frac{1}{2} \|u_n\|_\infty, \ \forall x\in  B_n,
\end{equation*}
and there exists a sequence $y_n\in\p B_n$ 
such that
\begin{equation}\label{yn}
u_n(y_n)= \frac{1}{2} \|u_n\|_\infty. 
\end{equation}
Choosing $x_n$ and $y_n$ in \eqref{xn}, \eqref{yn} respectively, $|u_n(x_n)-u_n(y_n)|=\frac12 \|u_n\|_\infty$ and using \eqref{Morrey}
\begin{equation*}
\frac12 \|u_n\|_\infty \le C   \|u_n\|_{W^{1,q^*}(\Om)} \, R_n^{2-N/q},
\end{equation*}
so
$$
R_n^{2-N/q}\geq Ch\big(\|u_n\|_\infty\big)^{-1/q} \|u_n\|_\infty.
$$
consequently
$$
R_n^{N}\ge C
\left( h\big(\|u_n\|_\infty\big)^{-1/q} \|u_n\|_\infty \right)^{\frac1{\frac{2}{N}-\frac{1}{q}}}.
$$
Observe that for $s\ge s_1,$

\begin{align}
h(s)^{-1/q} s
&= s^{-[(2^*-2)-2^*/q]}L(s)^{-1}\ \big(s|L'(s)|\big)^{1/q}\nonumber \\
&= s^{-2^*[2/N-1/q]}L(s)^{-1}\ \big(s|L'(s)|\big)^{1/q},\nonumber 
\end{align}
hence
\begin{align*}
R_n^N\ge  \|u_n\|_{\infty}^{-2^*}\, \Big(\big[ \|u_n\|_\infty\ L'_n\big]^{1/q}\,L_n^{-1} \Big)^{1/[2/N-1/q]},  
\end{align*}
where we have  denoted by
\begin{equation*}
 L_n:=L\big(\|u_n\|_\infty\big) ,\qquad
 L'_n:=\big|L'\big(\|u_n\|_\infty\big)\big|.
\end{equation*}

On the other hand, since $s\to s^{2^*+1}|L'(s)|$ is $RV_{2^*}$, by \eqref{bef5} in Corollary \ref{cor:RV:db}, there exists $C>0$ 
$$ u_n^{2^*+1}(x)\big|L'(u_n(x))\big|\geq C\|u_n\|_\infty^{2^*+1} L'_n,\quad \forall x\in B_n$$
so
\begin{align*}
C&\ge \int_\Om u_n^{2^*+1}|L'(u_n)|\, \,dx
\ge \int_{B_n}u_n^{2^*+1}|L'(u_n)|\, dx \\[.2cm]
&\ge C\,  \|u_n\|_\infty^{2^*+1}L'_n\, R_n^N\\[.2cm]
&\ge C  \, \|u_n\|_{\infty}\, L'_n \, \Big( \big( \|u_n\|_\infty\ L'_n \big)^{1/q}\,L_n^{-1} \Big)^{1/[2/N-1/q]}\\
&=  C \Big(\big( \|u_n\|_\infty\ L'_n \big)^{2/N}\,L_n^{-1}\Big)^{1/[2/N-1/q]}
\end{align*}
Consequently, 
$$
\dfrac{ \|u_n\|_\infty\ L'_n}{L_n^{\frac{N}{2}}}\le C ,
$$ 
contradicting the hypothesis \eqref{hypothesis1}. 
\end{proof}

\section{Proof of Theorem \ref{th:apriori}}
\label{sec:5}
The conclusions of Step 1 and Step 2 of the proof of Theorem \ref{th:apriori:cnys}, say \eqref{bdd:3} and \eqref{F:f:4}, remain valid.
 \smallskip
 
From Theorem \ref{usfF}.(i) for $q<2^*-1$ (see \eqref{eq:usfF}), there exists a constant $s_0>0$  such that
\begin{equation*}
2NF(s)-(N-2)sf(s)\ge \frac{N-2}2 \left(\frac{2^*}{q+1}-1\right)\,  sf(s),\ \forall  s>s_0.
\end{equation*}
Applying this inequality to any positive solution, and integrating on $\Om$ we obtain that
\begin{equation*}
2N\int_\Om F(u)\, \,dx -(N-2)\int_\Om uf(u)\, \,dx\ge 
C_0\int_\Om uf(u)\, \,dx-C_1,
\end{equation*}
for some constants $C_0,\ C_1$ independent on $u$. 

This,  and \eqref{F:f:4} yield
\begin{equation*}
\int_\Om u f(u)\, \,dx \le C,
\end{equation*}
for some constant $C$ independent of $u$. Consequently,
\begin{equation*}
\int_\Om |\na u|^2\, \,dx =	\int_\Om u f(u)\, \,dx \le C,
\end{equation*}
where $C$ is independent of $u$.
So, the $H^1_0(\Om)$-norm of any weak solution to \eqref{eq:ell:pb} is uniformly bounded. Now, \cite[Theorem 1.5]{Pardo_JFPTA_2023} ends the proof.

\begin{appendix}

\section{Examples of slowly  varying functions at infinity}
\label{app:A}
\renewcommand{\thesection}{\Alph{section}}
\numberwithin{equation}{section}

Slowly  varying functions are stable by multiplication, as well as by raising to any real power, see Proposition \ref{pro:L:al}. Typical examples of functions $L$ with slow variation at  $\infty$, and having unbounded or oscillating behaviors at $ \infty$, or both, are contained in the  Example \eqref{example}.

In this Appendix, we will check that the functions $L_i$ for $i=1,\cdots ,14$, 
defined in the following Example \eqref{example}, are slowly  varying functions at infinity, in other words, $L_i\in RV_0$.\\

We will say that the  functions  $g_1,\ g_2$ are {\it asymptotically similar at infinity,} and we will denote it by $g_1 \simeq_{+\infty} g_2$, when
$$
0 <\liminf_{s\to+\infty}\frac{g_1(s)}{g_2(s)}\le \limsup_{s\to+\infty}\frac{g_1(s)}{g_2(s)}<+\infty.
$$

\begin{example}\rm
\label{example}
The following set of functions satisfy $L_i\in RV_0$.

\begin{enumerate}
\item\label{1}  
$
L_{1}(s):=\left[\log (K + s)\right]^\al ,\ $  for $\al\ne 0,\ K > 1 $;

\item\label{2}
$
L_{2}(s):=\Big[ \log (K+s)/\log\big(K +\log (K+s)\big)\Big]^\al,\ $  for $\al\ne 0,\ K > 1 $;

\item\label{3}  
$ 
L_{3}(s):=\Big(\log_m (K + s)\Big)^\al ,
$  
where $\log_m$ represents the iterated logarithms for $m=1,2,3,\cdots$, defined in the following way. Let $\ \log_{1} (K + s)=\log (K + s),$ and
$$ 
\log_m (K + s)=\log\big(K +\log_{m-1}(K+s)\big)\q{for} \al\ne 0,\ K>1 ;
$$

\item\label{4} 
$L_{4}(s):=\exp{\Big(\al\,\log^\be (K+s)\Big)}$, for $\al\ne 0,\ 0<\be<1,\ K > 1$;  

\item\label{5} 
$
L_{5}(s):=\exp \Big(\al\, \log (K+s)/\log\big(K +\log (K+s)\big)\Big),\q{for} \al\ne 0,\ K > 1;
$

\item\label{6}
$
L_{6}(s):=\exp{\Big(\al\,\log^\be \big(\log (K+s)\big)\Big)},\quad
$
for $\al\ne 0,\ \be >0 $;  

\item\label{7}  
$
L_{7}(s):=\exp{\bigg(\al\,\D\prod_{i=1}^m \log_i ^{\be_i}(K_i+s)\bigg)}
$, 
for $\al\ne 0,\ 0<\be_i<1,\ K_i > 1,\ m=1,2,3,\cdots$;

\item\label{8}
$
L_{8}(s):=\exp \Big(\al\, \log (K+s)/(K+s)^\be \Big),$  for $ \al\ne 0, \be >0, K > 1;
$

\item\label{9}
$
L_{9}(s):=1 + \g \cos\Big[\log^\be (K+s)\Big],$  for $0<\be<1, K > 1,  |\g| < 1;
$

\item\label{10}
$
L_{10}(s):=1 + \g \cos\Big[\log \big (\log (K+s)\big)\Big],$  for $ K > 1,  |\g| < 1;
$

\item\label{11}  
$
L_{11}(s):=\exp \bigg[\al  \Big( \g\cos \big(\log^\be (K+s)\big)+ \log^\be (K+s)\Big) \bigg],
$ 
with  $\al\ne 0,\   |\g| <1,\ 0<\be<1,\ K>1$;

\item\label{12}  
$
L_{12}(s):=\exp \bigg[\al\,  \bigg(\g\cos \Big(\log \big (\log (K+s)\big)\Big)+ \log \big (\log (K+s)\big)\bigg)\bigg],
$
with  $\al\ne 0,\  |\g| < 1,\  K>1$;

\item\label{13}  
$
L_{13}(s):=\exp \bigg[\al  \bigg(1 + \g\cos \big(\log^\be (K+s)\big)+ \log^\be (K+s)\Big)\, \times\, \log^\be (K+s) \bigg],$ with  $\al\ne 0,\   |\g| < 1,\ 0<\be<1/2,\ K>1$;

\item\label{14}  
$
L_{14}(s):=\exp \bigg[\al\,  \bigg(1 + \g\sin \Big(\log \big (\log (K+s)\big)\Big)+ \log \big (\log (K+s)\big)\bigg)\,\times\, \log \big (\log (K+s)\big)\bigg],
$ 
with  $\al\ne 0,\  |\g| < 1,\  K>1$.
\end{enumerate}
\end{example}

\begin{proof}
Fixing $\tau>0$, we denote by $I=I_\tau :=\big[\min\{1,\tau\},\max\{1,\tau\}\big]$.

\noindent (1) We check it for the function
\begin{equation}\label{L1^}
\hat{L}_{1}(s):=L_{1}(s-K),    
\end{equation}
and observe that $\frac{\hat{L}_{1}(s)}{L_{1}(s)}\to1$ as $s\to+\infty$. 
For any $\tau, s>0$ we have 
\begin{align}\label{log:al}
	&\log^\al (\tau s)=\big(\log \tau + \log s\big)^\al=
	\log^\al(s)\,\big(1+ \log \tau /\log s\big)^\al
\end{align}
Then, 
$
\lim_{s\to +\infty}\frac{L_{1}(\tau s)}{L_{1}(s)}=1.
$
Besides, changing the variable $t=1/\log s$, and using l'Hôpital  it can be checked that, whenever $\al<1$
\begin{align}\label{lim:L1}
	&\lim_{s\to+\infty}\dfrac{\log^\al(s)\,\Big[\big(1+ \log \tau /\log s\big)^\al-1\Big]}{\log \tau}
	=\lim_{t\to 0^+}\frac{\big(1+ (\log \tau) t\big)^\al-1}{\log \tau\, t^\al}\\
	&=\lim_{t\to 0^+}\frac{ \al \big(1+ (\log \tau) t\big)^{\al-1}}{\al t^{\al-1}}= \lim_{t\to 0}\frac{t^{1-\al}}{ \big(1+ (\log \tau) t\big)^{1-\al}}\nonumber
	=0.
\end{align}
Hence, $\lim_{s\to+\infty}  \hat{L}_{1} (\tau s)-\hat{L}_{1} (s)=0$ and also
\begin{align}\label{dif:L1}
&\lim_{s\to+\infty}L_1 (\tau s)-L_1 (s)=0\qq{whenever} \al<1.
\end{align} 
\medskip

\noindent (2)  Example \ref{1} and Proposition \ref{pro:L:al}, parts (iii) and (iv) prove that $L_2\in RV_0$.
Moreover, it can be checked that
\begin{align}\label{dif:L2}
	&\lim_{s\to+\infty}L_2 (\tau s)-L_2 (s)	=0.
\end{align} 
\medskip

\noindent (3) Example \ref{1}  and Proposition \ref{pro:L:al}.(iv) proves the result.
It can also be checked that
\begin{align}\label{dif:L3}
\lim_{s\to+\infty}  L_{3} (\tau s)-L_{3} (s)=0.
\end{align} 
From definition, we also observe that
\begin{align}\label{L'7}
L_{3}' (s)&=\frac{\al}{K+s}\,
\frac{\log_m (K+s)^{\al-1}}{\D\prod_{k=1}^{m-1}\big(K+\log_k(K+s)\big)}\simeq_\infty \frac{\al}{s}\,\frac{L_{3}(s)}{\D\prod_{k=1}^{m}\log_k(s)}.
\end{align}
\medskip

\noindent (4)
Substituting  in \eqref{log:al}, \eqref{dif:L1} the exponent $\al$ by $\be,$ and using Example \eqref{1},   we get that $\lim_{s\to+\infty}\log L_4 (\tau s)-\log L_4 (s)=0$ whenever $\be<1.$ Now, Proposition \ref{pro:L:al}.(v) proves the result.
\medskip

\noindent (5) Example \eqref{2}, equation \eqref{dif:L2}, and Proposition \ref{pro:L:al}.(v) proves the result.
\medskip

\noindent (6)  Let   $L_{6}(s)$ be defined by \eqref{6}. Let  
$
\hat{L}_{6}(s):=s^{\al  [\log^\be  (\log  s)]/ \log  s}=\exp{\Big(\al\,\log^\be (\log (s))\Big)},\quad
$
for $\al\ne 0,\ \be >0 $.  
Taking logarithms, it can be seen that 
$
\lim_{s\to +\infty} \frac{\hat{L}_{6}(s)}{L_{6}(s)}=1.
$
Substituting  in  \eqref{dif:L1}  $\al$ by $\be,$ 
choosing $m=2$ in the definition of $ \hat{L}_3$, and using  \eqref{dif:L3}
we obtain that $\lim_{s\to+\infty}\log L_6 (\tau s)-\log L_6 (s)=0.$ Proposition \ref{pro:L:al}.(v) proves that $L_6\in RV_0$.
\medskip

\noindent (7) 
Let us check that  $\rho_{7}(s):=\frac{sL_{7}' (s)}{L_{7}(s)}\to 0$ as $s\to+\infty$, and then \cite[Theorem 2]{Lamperti_1958a} (see also Proposition \ref{pro:lamp}), will imply that $L_{7}\in RV_0$. Indeed, using \eqref{L'7}, we can write
$$
L'_{7}(s)=\al  
\sum_{i=1}^m \frac{\be_i}{K_i+s}\, \frac{\D\left(\prod_{j=1}^m\log_j^{\beta_j} (K_j+s)\right) \, L_{7}(s)}{\log_i (K_i+s)\D\prod_{k=1}^{i-1}\big(K_i+\log_k(K_i+s)\big)}
$$
Moreover, by definition, for all $\g>0$, and all $\Tilde{K}_{i},\ \Tilde{K}_{j}>1$,
$$
\frac{\log_m (\Tilde{K}_{i} + s)}{\log_{m-1}^\g(\Tilde{K}_{j}+s)}=\frac{\log\big(\Tilde{K}_{i} +\log_{m-1}(\Tilde{K}_{i}+s)\big)}{\log_{m-1}^\g(\Tilde{K}_{j}+s)}\to 0 \hbox{  as  } s\to+\infty,
$$
consequently, for any $k<j$, all $\Tilde{K}_{i},\ \Tilde{K}_{j}>1$, and all $\g>0$
\begin{equation*}
	\frac{\log_j (\Tilde{K}_{i} + s)}{\log_{k}^\g(\Tilde{K}_{j}+s)}\to 0 \hbox{  as  } s\to+\infty,    
\end{equation*}
and since $0<\be_j<1,$
\begin{equation*}
	\frac{\D\prod_{j=1}^m\log_j^{\beta_j} (K_j+s)}{\D\prod_{k=1}^{i}\log_k(K_i+s)} \to 0 \hbox{ as } s\to+\infty.
\end{equation*}
Hence,  since $0<\beta_i<1$ for all $1\leq i\leq m$,
\begin{align}
	|\rho_{7}(s)|&  \simeq_{+\infty}|\al |\sum_{i=1}^m \be_i\, \frac{\D\prod_{j=1}^m\log_j^{\beta_j} (K_j+s)}{\D\prod_{k=1}^{i}\big(K_i+\log_k(K_i+s)\big)}\to 0 \hbox{ as } s\to+\infty.\nonumber
\end{align}

\noindent (8)   
It is easy to get that for $
\al\ne 0,\ \be >0,\ K > 1
$
\begin{equation}\label{lim:inf:L9}
\lim_{s\to +\infty} \exp \Big(\al\, \log (K+s)/(K+s)^\be\Big)=1,
\end{equation}
that
$
\ L_{8}(s)=(K+s)^{\al/(K+s)^\be},$ and that $\lim_{s\to +\infty} \frac{\hat{L}_{8}(s)}{L_{8}(s)}=1,$ where $\hat{L}_{8}(s):=s^{\al/s^\be}
$
for $s>1$.
Moreover,
\begin{align*}
&\lim_{s\to +\infty}\frac{\hat{L}_{8}(\tau s)}{\hat{L}_{8}(s)}
= \lim_{s\to +\infty} \frac{(\tau s)^{1/(\tau s)^\be}}{s^{1/s^\be}}
= \lim_{s\to +\infty}\left( \frac{\tau s}{s}\right)^{1/(\tau s)^\be}\,
s^{1/s^\be\, ( 1/\tau ^\be-1)}
\nonumber\\
&= \lim_{s\to +\infty}\tau ^{1/(\tau s)^\be}\,
\,
s^{1/s^\be\, ( 1/\tau ^\be-1)}=\tau^0 1^{1/\tau ^\be-1}=1.
\end{align*} 
\medskip

\noindent (9)
Denoting $\hat{L}_9 (s)= L(s-K)=1 + \g \cos\big[\log^\be (s)\big],$ we check that $\hat{L}_9\in RV_0$.  Observe that $\hat{L}_9$ is positive since $|\gamma|<1$.
Substituting  in \eqref{log:al}  $\al$ by $\be,$  we can write
\begin{align*}
\cos\big(\log^\be (\tau s)\big)
&=\cos\big(\log^\be(s)\big)\, \cos\bigg(\log^\be(s)\,\Big[\big(1+ \log \tau /\log s\big)^\be-1\Big]\bigg)\\
&\quad - \sin\big(\log^\be(s)\big)\, \sin\bigg(\log^\be(s)\,\Big[\big(1+ \log \tau /\log s\big)^\be-1\Big]\bigg)\nonumber
\end{align*}
Now,  \eqref{lim:L1} gives that
\begin{equation}\label{dif:L9}
\lim_{s\to+\infty}\frac{\hat{L}_9 (\tau s)}{\hat{L}_9 (s)}=1 \quad \hbox{ and }  \lim_{s\to+\infty}\hat{L}_9 (\tau s)-\hat{L}_9 (s)=0.
\end{equation}
\medskip

\noindent (10)
Consider $
\hat{L}_{10}(s):=L_{10}(s-K)=1 + \g \cos\Big[\log \big (\log (s)\big)\Big],$  for $  |\g| < 1,
$ and observe that $\frac{\hat{L}_{10}(s)}{L_{10}(s)}\to1$ as $s\to+\infty$, 
and let us check that 
$\hat{L}_{10}$ is  $RV_0$ . Indeed, for $|\gamma|<1$ this function is positive and for any $\tau>0$ and $s>0$ we have
\begin{align*}
\log\big(\log (\tau s)\big)=\log\log s + \log\big(1+ \log \tau /\log s)\big),
\end{align*}
and $\log\big(1+ \log \tau /\log s\big)\to 0$, so, as $s\to+\infty,$ 
\begin{equation*}
\cos\big(\log\big(1+ \log \tau /\log s)\big)\to 1, \quad \sin\big(\log\big(1+ \log \tau /\log s)\big)\to 0
\end{equation*}
Hence $\lim_{s\to+\infty}\frac{\hat{L}_{10} (\tau s)}{\hat{L}_{10} (s)}=1 $ and $\lim_{s\to+\infty}\hat{L}_{10} (\tau s)-\hat{L}_{10} (s)=0,$ and so
\begin{equation}\label{dif:L10}
\lim_{s\to+\infty}\frac{L_{10} (\tau s)}{L_{10} (s)}=1 \quad \hbox{ and }  \lim_{s\to+\infty}L_{10} (\tau s)-L_{10} (s)=0.    
\end{equation}

\medskip

\noindent (11)
When   $L_{11}(s)$  is given by \eqref{11}, we  compute the corresponding limit, 
for $\hat{L}_{11}(s)=L_{11}(s-K)$.
Taking logarithms, and using \eqref{dif:L1} and \eqref{dif:L9}
\begin{align*}
& \lim_{s\to+\infty}\, \frac{1}{\al}\, \log \left(\frac{\hat{L}_{11} (\tau s)}{\hat{L}_{11} ( s)}\right)\\
&=\lim_{s\to+\infty}\gamma \Big(\cos (\log^\beta (\tau s)) -\cos (\log^\beta (s)\Big) + \log^\beta (\tau s) -\log^\beta (s)\\
&=\lim_{s\to+\infty} \gamma (\hat{L}_{9} (\tau s)-\hat{L}_{9} (s) )+ \hat{L}_{1} (\tau s)-\hat{L}_{1} (s)=0,
\end{align*}
and then $L_{11}\in RV_0.$ 
\medskip

\noindent (12) 
For  $L_{12}(s)$   given by \eqref{12}, we  compute the corresponding limit, as in \eqref{L1^}, for $\hat{L}_{12}(s):=L_{12}(s-K)$.
Taking logarithms,  for any $\tau>0$,
\begin{align*}
&\lim_{s\to +\infty}\left|\log \left(\frac{\hat{L}_{12}(\tau s)}{\hat{L}_{12}(s)}\right)\right|
=\lim_{s\to +\infty}|\al| \bigg|\Big( \g\cos\big(\log\log (\tau s)\big)+\log\log (\tau s)\Big)
\\
&\qquad\quad -\Big( \g\cos \big(\log\log s\big) + \log\log s\Big) \bigg|\\
&\quad =\lim_{s\to+\infty} |\alpha|\big( \hat{L}_{10} (\tau s)-\hat{L}_{10} (s) + \hat{L}_{3} (\tau s)-\hat{L}_{3} (s)\big)=0,
\end{align*} 
choosing $m=2$ in the definition of $ \hat{L}_3$, and using \eqref{dif:L10} and \eqref{dif:L3}.
\medskip

\noindent (13) When   $L_{13}(s)$  is given by \eqref{13}, we will compute the corresponding limit, as in \eqref{L1^}, for $\hat{L}_{13}(s)=L_{13}(s-K)$.
Taking logarithms,  and  using the mean value theorem,
for any $\tau>0,\ |\g|<1$,  and  $0<\be<1/2$
\begin{align*}
&\lim_{s\to +\infty}\left|\log \left(\frac{\hat{L}_{13}(\tau s)}{\hat{L}_{13}(s)}\right)\right|
\nonumber\\
&\ =\lim_{s\to +\infty} |\al|\Big|\big[1 + \g\cos \big(\log^\be (\tau 
s)\big)+\log^\be (\tau 
s)\big]\,\log^\be (\tau s) \\
&\qquad\qquad -\big[1 + \g\cos \big(\log^\be (s)\big)+\log^\be (s)\big]\,\log^\be s \Big|\nonumber\\
&\ \le\be\,\lim_{s\to +\infty}\max_{\te\in I} \left\{\dfrac{\Big(2\big(1+|\g|\big)+\log^{\be} (\te s)\Big)\log^{\be-1} (\te s)}{\te }\right\}\,|\tau-1|
=0.
\end{align*} 
\medskip

\noindent (14) 
For  $L_{14}(s)$   given by \eqref{14}, we  compute the corresponding limit, as in \eqref{L1^}, for $\hat{L}_{14}(s):=L_{14}(s-K)$.
Taking logarithms,  and  using the mean value theorem,
for any $\tau>0$,  
\begin{align*}
&\lim_{s\to +\infty}\left|\log \left(\frac{\hat{L}_{14}(\tau s)}{\hat{L}_{14}(s)}\right)\right|\nonumber\\
&=\lim_{s\to +\infty}|\al| \bigg|\Big(1 + \g\sin\big(\log\log (\tau s)\big)+\log\log (\tau s)\Big)
\, \log\log (\tau s)
\\
&\quad \quad -\Big(1 + \g\sin \big(\log\log s\big) + \log\log s\Big) \log\log (s)\bigg|\\
&\ \le|\al| \,\lim_{s\to +\infty}\max_{\te\in I} \left\{\frac{(|\g|+1)\log\log (\te s)+\Big(1+|\g|+\log\log (\te s) \Big)}{\te s\log (\te s)}\right\}\\
&\quad \quad \,\times\,|\tau-1|=0.
\end{align*} 
\end{proof}

\section{Examples of \texorpdfstring{$L^\infty(\Om)$}  {\it a priori} bounds for regularly varying functions at infinity}
\label{app:B}

Let us denote 
\begin{align*}
&f_i(s)=s^{2^*-1}L_i(s), \qq{with} s> 0,\quad \al<0,
\end{align*}
for  $L_i\in RV_0$,  for $i=1,\cdots ,14,$ depending on $\al,$ and possibly more parameters, described in the Example \ref{example}.  

In this Appendix, we will check if the positive solutions to \eqref{eq:ell:pb} for the nonlinearities associated to the previously defined slowly varying functions $L_i$ for $i=1,\cdots ,14$, 
are $L^\infty(\Om)$  {\it a priori} bounded.\\
We use Theorem \ref{th:apriori:cnys},  specifically that  
$f_i$ satisfying {\rm \ref{f1}}-\ref{f3}, 
joint with
\begin{equation}\label{L:inf:h:2}
\lim_{s\to+\infty} \,\frac{s\,L_i'(s)}{L_i^{N/2}(s)}=+\infty\implies \|u\|_\infty\le C.
\end{equation}

\begin{proof}[Proof of Corollary \ref{cor:ex}]
\noindent (1)  Hypothesis \ref{f1}--\ref{f3} and \eqref{L:inf:h:2}  are satisfied for $L_{1}(s)$  given by \eqref{1}, whenever $|\al|>\frac{2}{N-2}.$ Indeed, 
$
L'_{1}(s)=\alpha \,\dfrac{{\log^{\alpha -1} \left(K+s\right)}}{K+s},
$
and so, 
$$
\lim_{s\to+\infty} \,\frac{L_{1}^{N/2}(s)}{sL'_{1}(s)}=
\lim_{s\to+\infty} \, \frac{1}{\alpha}\left[\log (K + s)\right]^{1-\frac{|\al|(N-2)}{2}} =0
\iff |\al|>\frac{2}{N-2}.\checkmark
$$
Moreover, $f_{1}$ satisfies \ref{f1}--\ref{f3}  for $\alpha<0$.
\medskip

\noindent (2)  Hypothesis \ref{f1}--\ref{f3} and \eqref{L:inf:h:2}  are satisfied for $L_{2}(s)$  given by \eqref{2}. Indeed, 
$ 
L'_{2}(s)	\simeq_{+\infty}\, \frac{\alpha}{s}\, L_{2}(s)\, \frac{1}{\log(s)},
$ 
and so, 
$$
\lim_{s\to+\infty} \,\frac{L_{2}^{N/2}(s)}{sL'_{2}(s)}=
\lim_{s\to+\infty} \, \dfrac{\log\big(K+\log(K+s)\big)^{\frac{|\al|(N-2)}{2}}}{\alpha\,\left[\log ( s)\right]^{\frac{|\al|(N-2)}{2}-1}}=0
\iff |\al|>\frac{2}{N-2}.
$$
Also, $f_{2}$ satisfies \ref{f1}--\ref{f3}  for $\alpha<0$.
\medskip

\noindent (3) 
for   $L_{3}(s) =\big(\log_m (K + s)\big)^\al$  defined by \eqref{3} with $m\ge 2$, hypothesis \eqref{hypothesis1}  is not satisfied for any $\al<0$.
Indeed, from definition, we obtain  
\begin{align*}
	L_{3}' (s)
	&\simeq_{+\infty}\al\,\frac{1}{K+s}\,\frac{L_{3}(s)}{ \D\prod_{k=1}^{m}\log_k(s)},\qq{and} |L_{3}'|=-L_{3}' \in RV_{-1}.
\end{align*}
Also
\begin{align*}
	\frac{L_{3}(s)^{N/2}}{s|L_{3}' (s)|}
	&\simeq_{+\infty}|\al|^{-1}\, \big(\log_m (K + s)\big)^{\al (N-2)/2}\, \D\prod_{k=1}^{m}\big(K+\log_k(K+s)\big) 
\end{align*}
and hence
$$
	\lim_{s\to+\infty} \,\frac{L_{3}^{(N-2)/2}(s)}{|\rho_{3}(s)|}= +\infty \qq{for any} \al\ne 0. 
$$ 
Consequently, hypothesis \ref{f1}-\ref{f3} are satisfied, but not \eqref{L:inf:h:2}.
\medskip

\noindent (4) For   $L_{4}(s)$  defined by \eqref{4}, 
$
L'_{4}(s) 
=\alpha \,\beta \,\dfrac{{\log^{\beta -1} \left(K+s\right)}}{K+s}\,L_{4}(s),
$ 
and so, for $\al<0$ and $0<\be<1$,
$$
\lim_{s\to+\infty} \,\frac{L_{4}^{N/2}(s)}{sL'_{4}(s)}
=\lim_{s\to+\infty} \,\frac{\exp{\Big(-\frac{|\al|(N-2)}{2}\log^\be (K+s)\Big)}}{ \alpha \beta\log^{\be-1} (K+s)}=0.\ \checkmark
$$
Moreover,  $f_{4}$ satisfies {\rm \ref{f1}}-\ref{f3} for $\al<0$ and $0<\be<1$.
\medskip

\noindent (5) If   $L_{5}(s)$ is defined by \eqref{5},  
$
	L'_{5}(s)
    \simeq_\infty\dfrac{\al}{s} \
	\dfrac{L_{5}(s)}{\log (s) \log\big(\log (s)\big)}\, ,
$ 
and so, 
$$
\lim_{s\to+\infty} \,\frac{L_{5}^{N/2}(s)}{sL'_{5}(s)}
=\lim_{s\to\infty} \, \dfrac{1}{\al}\,L_{5}^{(N-2)/2}(s)(s)\, \log (s) \log\big(\log (s)\big)=0. \checkmark 
$$
Moreover, $f_{5}$ satisfies \ref{f1}-\ref{f3} for $\alpha<0$. 
\medskip

\noindent (6) For   $L_{6}(s)$  defined by \eqref{6}, we have 
$ 
\hat{L}'_{6}(s) =\alpha \,\beta \,\dfrac{{\log^{\beta -1}\left(\log(s)\right)}\,}{s\,\log(s)}\, \hat{L}_{6}(s),
$ 
and so, 
$$
\lim_{s\to+\infty} \,\frac{\hat{L}_{6}^{N/2}(s)}{s\hat{L}'_{6}(s)}=
\lim_{s\to+\infty} \,\frac{\exp{\Big(\al(N-2) \,{\log^\beta \left(\log(s)\right)}/2\Big)}\, \log(s)}{ \alpha \beta{\log^{\beta -1}(\log s)}}=0, \hfill \checkmark
$$
Moreover, $f_{6}$ satisfies {\rm \ref{f1}}-\ref{f3} for $\al<0$,  and  $0<\beta<1$.
\medskip

\noindent (7) 
For   $L_{7}(s)$ defined by \eqref{7}, 
\begin{align*}
	L'_{7}(s)&=\al  
	\sum_{i=1}^m \frac{\be_i}{K_i+s}\, \frac{\D\left(\prod_{j=1}^m\log_j^{\beta_j} (K_j+s)\right) \, L_{7}(s)}{\log_i (K_i+s)\D\prod_{k=1}^{i-1}\big(K_i+\log_k(K_i+s)\big)} 
\end{align*}
and
\begin{align*}
	\frac{s|L'_{7}(s)|}{L_{7}(s)}
    &  \simeq_{+\infty}|\al |\sum_{i=1}^m \be_i\, \frac{\D\left(\prod_{j=1}^m\log_j^{\beta_j} (s)\right)  }{\D\prod_{k=1}^{i}\log_k(s)}=:\Tilde{l}_{7}(s),\nonumber
\end{align*} 
and so, for any $\al<0$,
\begin{align*}
	& \frac{L_{7}^{N/2} (s)}{s|L'_{7} (s)|}\simeq_{+\infty}\frac{L_{7}^{\frac{N-2}{2}} (s)}{\Tilde{l}_{7} (s)}\\
	&\le_{+\infty} C
	\exp{\Big(\al\frac{N-2}{2}\,\D\prod_{j=1}^m\log_j^{\beta_j} (s)\Big)}\, \prod_{k=1}^m \big(\log_k(s)\big)
	\to 0  \hfill  \checkmark   
\end{align*}
as $s\to+\infty.$ Moreover, $f_{7}$ satisfies {\rm \ref{f1}}-\ref{f3}.
\medskip

\noindent (8) 
If   $L_{8}(s)$ is defined by \eqref{8}, then \eqref{lim:inf:L9} implies that $f_{8}$ is not slightly subcritical.\medskip

\noindent (9) \& (10):
For   $L_{9}(s)$  defined by \eqref{9}, observe that $\liminf_{s\to\infty} L_{9}(s)=1-|\g|>0$, and so $f_{9}$ is not slightly subcritical. Likewise, for  $L_{10}(s)$  defined by \eqref{10}, $f_{10}$ is not slightly subcritical.\medskip

\noindent (11) For   $L_{11}(s)$  defined by \eqref{11}, $\al< 0,\   |\g| < 1,\ 0<\be<1,\ K>1$;
$$
L'_{11}(s)=\al\be L_{11}(s)\,\frac{\log^{\beta-1} (K+s)}{K+s} \,\Big(1-\gamma \,\sin\big(\log^\beta (K+s) \big)\Big) ,
$$
then  $L_{11}'<0$ if, say,  $\alpha<0$, $\beta>0 $ and $|\gamma|<1$. Moreover
$
	\dfrac{s|L'_{11}(s)|}{L_{11}(s)}\simeq_{+\infty} |\alpha |\beta \log^{\be-1} (K+s)
$ 
and so, for any $\al<0$,
\begin{align*}
& \frac{L_{11}^{N/2} (s)}{s|L'_{11} (s)|}\simeq_{+\infty}\, \frac{1}{|\alpha| \beta } L_{11}^{\frac{N-2}{2}} (s)\, \log^{1-\be} (K+s)
\to 0 .  \hfill \checkmark 
\end{align*} 
Moreover, for  $\alpha<0$, $|\gamma|<1$ and $0<\beta<1 $, $f_{11}$  satisfies \ref{f1}-
\ref{f3}.\\

\noindent (12) For   $L_{12}(s)$  defined by \eqref{12}, with $\alpha<0$,  and $|\gamma|<1$, 
$$
L'_{12}(s)=\al L_{12}(s)\,\frac{\frac{1}{(K+s)}}{\log (K+s)} \,\Big[1-\gamma \,\sin\Big(\log \big (\log (K+s)\big)\Big)\Big]<0 .
$$
Moreover
$
\dfrac{s|L'_{12}(s)|}{L_{12}(s)}\simeq_{+\infty} |\alpha | \log^{-1} (K+s)
$
and so, for any $\al<0$,
\begin{align*}
& \frac{L_{12}^{N/2} (s)}{s|L'_{12} (s)|}\simeq_{+\infty}\, \frac{1}{|\alpha|} L_{12}^{\frac{N-2}{2}} (s)\, \log (K+s)
\to 0 ,  \text{ hence \eqref{L:inf:h:2} holds.} 
\end{align*} 
Moreover, $f_{12}$ satisfies \ref{f1}-\ref{f3}.\\

\noindent (13) Let   $L_{13}(s)$  be given by \eqref{13}, we act as in \eqref{1}, for
$
\hat{L}_{13}(s):=L_{13}(s-K)=\exp \bigg[\al  \Big(1 + \g\cos \big(\log^\be (s)\big) + \log^\be (s)\Big)\, \times\, \log^\be (s) \bigg],$ with  $0<\be<1/2,\ |\g|<1.
$
$
\hat{L}'_{13}(s) 
	 \simeq_{+\infty} \alpha \,\beta \,\hat{L}_{13}(s)\,\dfrac{\log^{2\beta-1} (s)}{s}
$ 
and  $|\hat{L}'_{13}|\in RV_{-1}$.  

On the other hand,
\begin{align*}
	\lim_{s\to\infty}\frac{\hat{L}_{13}(s)^{N/2}}{s|\hat{L}'_{13}(s)|}= \lim_{s\to\infty}\frac{\hat{L}_{13}(s)^{(N-2)/2}}{|\alpha| \beta\log^{2\beta-1} (s)}
	=0,\text{ so \eqref{L:inf:h:2} holds.}
\end{align*}
Moreover, $f_{13}$ satisfies {\rm \ref{f1}}-\ref{f3}.\\

\noindent (14)
Let us consider   $\hat{L}_{14}(s)=L_{14}(s-K)$, with $L_{14}$ defined by \eqref{14}, for $K>1,\  |\g| < 1,\ \be<1,\ \al<0,$
$
\hat{L}'_{14}(s)  \simeq_{+\infty}\hat{L}_{14}(s)\, 
\dfrac{\alpha\, \log(\log s)}{s\, \log (s)} ,
$ 
and so $|\hat{L}'_{14}|\in RV_{-1}$. 
Moreover
\begin{align*}
	\lim_{s\to+\infty}\frac{\hat{L}_{14}(s)^{N/2}}{s|\hat{L}'_{14}(s)|}=
	\lim_{s\to+\infty}\frac{\hat{L}_{14}(s)^{(N-2)/2}\, \log (s)}{ |\alpha|\log(\log s)}
	=0,\text{ and \eqref{L:inf:h:2} holds.}
\end{align*}
Also, $f_{14}$ satisfies {\rm \ref{f1}}-\ref{f3}  for any $\al<0$.
\end{proof}

\end{appendix}

\section{Acknowledgments} 


The second author thanks Professor Alfonso Castro for having introduced me to the topic of {\it a priori} bounds for elliptic equations, as well as for a vastitude of wonderful discussions. Collaborating with him has been a genuinely challenging and stimulating experience, and a pleasant way to move forward in the profession, I can not be grateful enough.


\end{document}